\documentstyle[12pt]{article}

\newcommand{\be}{\begin{equation}}
\newcommand{\ee}{\end{equation}}
\newcommand{\ba}{\begin{eqnarray}}
\newcommand{\ea}{\end{eqnarray}}
\newcommand{\baa}{\begin{eqnarray*}}
\newcommand{\eaa}{\end{eqnarray*}}
\newcommand{\bb}{\begin{thebibliography}}
\newcommand{\eb}{\end{thebibliography}}
\newcommand{\ci}[1]{\cite{#1}}

\newcommand{\lab}[1]{\label{#1}}
\newcommand{\re}[1]{(\ref{#1})}

\newcommand\fac[2]{\mbox{$\frac{#1}{#2}$}}

\newcommand\lbr{\lbrace}
\newcommand\rbr{\rbrace}

\begin{document}
\bibliographystyle{unsrt}
\begin{titlepage}

{\hfill November 2020}

\vspace*{25mm}

\begin{center}

{\LARGE \bf Hidden Symmetry of the Hahn\\[2mm]

Problem for the Quantum Algebra $sl_q(2)$}\\[5mm]

A.N.Lavrenov

\vspace{2mm}
{\it Physics and Mathematical Department, Belarus State Pedagogical University, Minsk, 220125 Belarus}\\[3mm]

(Submitted 01 November 2020)\\[5mm]
\end{center}

\vspace{5mm}

\begin{abstract}

A special case of Askey-Wilson algebra $AW(3)$ with three generators is shown to serve
as a hidden symmetry algebra underlying the Hahn problem for the quantum algebra $sl_q(2)$. On the base of this
hidden symmetry 
the  corresponding  Clebsch-Gordan coefficients 
in terms of the q-Hahn polynomials is found.

\medskip
PACS numbers: 02.30.Gp, 03.65Fd, 11.30.Na
\end{abstract}
\end{titlepage}

\newpage

\section{Introduction}

\vspace{8mm}

As is well known, any symmetry of the problem under consideration corresponds to its own 
symmetry operator. If there are several of them, then they can form 
an dynamic or hidden symmetry algebra. In particular, the Askey-Wilson algebra is considered 
as the most general algebra for problems with the  polynomial solutions \ci{zh1}.
Her special case of Hahn algebra has lately become quite popular. Even the  so-called 
meta-Hahn
algebra  has recently been introduced for a unified algebraic underpinning of the Hahn
polynomials and rational functions \ci{zh2}.
On the other hand, there is the problem of finding the algebra of dynamic or hidden symmetry 
from the general principles of its construction.
One of the approaches to solving this problem was proposed in \ci{zh3} for general systems possessing 
 $SU(1,1)\oplus SU(1,1)$  dynamical symmetry. 
The quadratic Hahn algebra $QH_(3)$ was shown to serve as a hidden symmetry in both quantum and 
classical pictures.        
Attempt of its q-generalization can be found in \ci{qH}, but they were all limited only to the 
$SU_q(1,1)$ case. Although it should be noted that
in \ci{zh4} a new addition rule is proposed for nonlinear algebras including $sl_q(2)\oplus sl_q(2)$  and 
two types of q-oscillator algebra.

The purpose of this paper is to present an analogous to \ci{zh3} algebraic treatment 
of hidden symmetry for all types of algebras obeying that 
a new addition rule which  was proposed in \ci{zh4}.

The paper is organized as follows. In Sec.II, we recall the addition rule
for different types of $sl_q(2)$ algebras in accordance with Ref. \ci{zh4}.
In Sec.III, the special case of Askey-Wilson algebra $AW(3)$ is shown to be the hidden symmetry algebra for
this case. 
On the base of this
hidden symmetry 
the  corresponding  Clebsch-Gordan coefficients 
in terms of the q-Hahn polynomials will be presented in Section IV.
Concluding remarks and perspectives will form the last section.

\section{Different types of $sl_q(2)$ and their addition rule}

This section provides the necessary background material on  the addition rule  for nonlinear algebras. 
In particular, 
we represent the next own notation for the $sl_q(2)$ algebra (compare with \ci{zh4, zh5, maj}), which is generated by three operators 
$A_0$, $A_{+}$, $A_{-}$ obeying the relations:
\ba
\lab{DefiningRelations}
\; [A_0,A_\pm]&=&\pm A_\pm,
\nonumber \\
\; [A_-,A_+]&=&g(A_0+1/2) - g(A_0- 1/2)=
\nonumber \\
\;                &=&(q-q^{-1})(a_1 q^{2 A_0}- a_2 q^{-2 A_0}),
\lab{def} \ea
where $[a,b]=ab-ba$; $g(x)=a_1 q^{2x}+a_2 q^{-2x}$

In what follows we shall denote the algebra $sl_q(2)$ with
commutation relations \re{def}  by the symbol $(a_2, a_1)$.

The special cases of $(a_2, a_1)$ algebra are:

(i)   $su_q(2)$ if      $a_2 = a_1>0$ and $q>1$ or $a_2 = a_1<0$ and $0<q<1$ ; 

(ii)   $su_q(1,1)$ if   $a_2 = a_1<0 $ and $q>1$ or $a_2 = a_1>0$ and $0<q<1$ ; 

(iii)  $cu_q(2)$ if    $a_2=-a_1<0$  and $q>1$ or $a_2 = -a_1>0$ and $0<q<1$ ; 

(iv)  $eu_q^+$ if      $a_2<0, a_1=0$ and $q>1$ or $a_2>0, a_1=0$ and $0<q<1$ ; 

(v)   $eu_q^-$ if       $a_2=0, a_1>0$ and $q>1$ or $a_2 =0,  a_1<0$ and $0<q<1$. 

The Casimir operator of the $(a_2, a_1)$ algebra which commutes with all generators has the expression
\ba
\hat Q    &=A_+A_-   -   g(A_0- 1/2) &= A_+A_-   -    a_1 q^{2A_0 - 1 } - a_2 q^{1 - 2A_0}                     =
\nonumber \\
\;    &=A_-A_+   -  g(A_0+1/2) &=A_-A_+   -  a_1 q^{2A_0+1} - a_2 q^{-2A_0-1}.
\lab{cas} \ea

In view of the defining relations (\ref{DefiningRelations}), it is clear that $sl_{q}(2)$ has a ladder representation. 
Let $\mu>0$ be a positive real number and consider the infinite-dimensional vector space  $V^{(\mu)}$ spanned by the orthonormal  basis vectors 
$e_{n}^{(\mu)}$, 
$n\in \left\{0, N\right\}$, and endowed with the actions
\be
\lab{Actions}
A_0 e_{n}^{(\mu)} = (n + \mu) e_{n}^{(\mu)},
\quad A_{+}e_{n}^{(\mu)} = r_{n+1} e_{n+1}^{(\mu)},
\qquad A_{-}e_{n}^{(\mu)} = r_{n} e_{n-1}^{(\mu)},
\ee
\ba
\lab{Actions}
A_0 e_{n}^{(\alpha)} = (n + \alpha) e_{n}^{(\alpha)},
\quad A_{+}e_{n}^{(\alpha)} = r_{n+1}^{(\alpha)} e_{n+1}^{(\alpha)},
\qquad A_{-}e_{n}^{(\alpha)} = r_{n}^{(\alpha)} e_{n-1}^{(\alpha)},
\ea
with $\langle e_i, e_j\rangle = \delta_{ij}$ and where $r_n$ is given by
\be
r_n^2=(q^{2n}-1)(a_1 q^{2\mu-1} - a_2 q^{1-2n-2\mu)})  = (q^{n}-q^{-n})(a_1 q^{2\mu-1+n}-a_2 q^{1-n-2\mu})
\lab{r}\ee
As expected from Schur's lemma, the Casimir operator $Q$ acts on canonical basis $e_{n}^{(\mu)}$ as a multiple of the identity:
\ba
\lab{CasimirAction}
\hat Q\,e_{n}^{(\mu)}=Q(\mu)e_{n}^{(\mu)} \equiv 
[- a_1 q^{2\mu - 1 } - a_2 q^{1 - 2\mu}] e_{n}^{(\mu)}.
\ea
Fixing the value of the Casimir operator $Q(\mu)$ we get a unitary representation of the
$(a_2, a_1)$ algebra. In this paper we restricted ourselves to the representations of the
positive discrete series $D^+_\mu$ where
$r_n>0$ and the state $e_{0}^{(\mu)}$ is the vacuum of
the representation $D^+_\mu$, i.e. $r_0=0$.

The $(a_2, a_1)$  algebra possesses an addition property that can be presented in the following way \ci{zh4}. Let $\{A_{0},A_{\pm},Q_{A}\}$ and $\{B_{0},B_{\pm},Q_{B}\}$ be two mutually commuting sets of $sl_{q}(2)$ generators and denote the corresponding algebras by $\mathcal{A}$ and $\mathcal{B}$. A third algebra, denoted $\mathcal{C}=\mathcal{A}\oplus\mathcal{B}$, is obtained by defining
\ba
\lab{res}
C_{0}&=&A_0+B_0,
\nonumber \\
C_\pm&=&A_\pm q^{-B_0}+B_\pm q^{A_0}.
\lab{add} \ea 

The addition rule \re{add} is the same as for ordinary $sl_q(2)$ algebra
\ci{maj}, however the algebras $(a_2, a_1)$ and $(b_2, b_1)$ in \re{add} may
have different types \ci{zh4,zh5}. It is easily seen that in order for the
operators $C_0, C_\pm$ to form new $(c_2, c_1)$ algebra,
the following relations must
be fulfilled:
\be
c_2=a_2, \qquad d \equiv a_1= b_2, \qquad c_1=b_1;\\
\lab{connect} \ee

In symbolic form the addition rule \re{add} can be written as
\be   
(a_2, d) \oplus (d, b_1) = (a_2, b_1)
\lab{symb} \ee

 The Casimir operator of the resulting algebra
\ba
\hat Q_{\mathcal{C}=\mathcal{A}\oplus\mathcal{B}}=C_+C_-   -   g(C_0- 1/2),
\ea
may be cast in the form
\ba
\lab{Cas-Full}
\hat Q_{\mathcal{C}=\mathcal{A}\oplus\mathcal{B}}=
\left\{
q A_+B_-  + q^{-1} B_+A_- +
(q +q^{-1}) d q^{A_0 - B_0 }  + \hat Q_B q^{2 A_0}  +  \hat Q_A q^{- 2 B_0}
\right\}
q^{A_0 - B_0 },
\ea
where $\hat Q_{i}$, $i\in \{A,B\}$, are the Casimir operators of the algebras $\mathcal{A}$ and $\mathcal{B}$.

\vspace{6mm}
\section{$AW(3)$ algebra and the Hahn problem}

According to the abstract schemes which was proposed in \ci{zh3} for general systems possessing   $SU(1,1)\oplus SU(1,1)$  dynamical symmetry for the Hamiltonian $C_0=A_0 +B_0$ we have two independent  integrals $K_1, K_2$ commuting with $C_0$: the difference between the original operators $K_1 =\Delta \equiv A_0 - B_0$ and Casimir
operators for the resulting algebra $K_2=\hat Q_{\mathcal{C}=\mathcal{A}\oplus\mathcal{B}}$.
The Hahn problem consists in finding the 
overlaps coefficients between the eigenbases of those two operators $K_1, K_2$ and
corresponds to the Clebsch-Gordan problem for $SU(1,1)$.
This problem is non-trivial because
the operators $K_1$ and $K_2$  do  not commute with one another.

Here we aim to construct a q-deformation of the above abstract scheme, preserving the general algebraic foundations for this approach, i.e. to solve  the Hahn problem for the quantum algebra $sl_q(2)$. 
In other words, for our case we  introduce 
one of them symmetry operator $K_2$ as  Casimir
operators for the resulting algebra  $(c_2, c_1)=(a_2, d) \oplus (d, b_1) = (a_2, b_1)$, i.e $K_2=\hat Q_{\mathcal{C}=\mathcal{A}\oplus\mathcal{B}}$:
\ba
K_2    &=C_+C_-   -   g(C_0- 1/2) &= (A_+ q^{-B_0}+B_+ q^{A_0}) 
(A_- q^{-B_0}+B_- q^{A_0}) -   c_1 q^{2C_0 - 1 } - c_2 q^{1 - 2C_0}=                     ,
\nonumber \\
\;    &=(\pm)q^{ \Delta } +P_2(q^{ \Delta }) &=  (q A_+B_-  + q^{-1} B_+A_-) q^{ \Delta } +
p_2 q^{2  \Delta }  + p_1 q^{ \Delta } + p_0,
\lab{13} \ea
where $p_2 = (q +q^{-1}) d;\quad p_1 = (Q_B q^{C_0}  + Q_A q^{- C_0}); \quad p_0=0$.

It is clear that the second symmetry operator $K_1$  is, in the general case, 
some function of  $\Delta$, i.e $K_1= f (\Delta)$.
The formula \re{13} suggests an explicit form of this function $f$ in the following form
\be
K_1    = q^{k \Delta},
\lab{14}
\ee
where while $k$ is arbitrary real parameteris.

Consider the following sum $K_1^2  K_2 + K_2 K_1^2$ in the two-operator approach to the Askey-Wilson algebra \ci{la1,daha}, trying to transform the each term of sum into the term $K_1  K_2 K_1$:
\ba
K_1^2  K_2 + K_2 K_1^2  &=&(q^{2k} +q^{-2k}) K_1  (q A_+B_-  + q^{-1} B_+A_-) q^{ \Delta } K_1  + 2P_2(q^{ \Delta }) K_1^2,
\nonumber \\
\;   &=&(q^{2k} +q^{-2k}) K_1 (K_2  - P_2) K_1  + 2P_2 K_1^2,
\nonumber \\
\;   &=&(q^{2k} +q^{-2k}) K_1 K_2  K_1  - (q^{k} -q^{-k})^2 P_2 K_1^2.
\lab{15} \ea

To remain in quadratic combinations of the original operators $K_1,  K_2$ for the last term $P_2 K_1^2$, it is necessary to take the value of $k$ equal to -1: 
\be
K_1^2  K_2 + K_2 K_1^2  - (q^{2} +q^{-2}) K_1 K_2 K_1+(q -q^{-1})^2 (p_2  + p_1  K_1)= 0,
\lab{16} \ee

Consider in a similar way the sum $K_2^2  K_1 + K_1 K_2^2$, trying to transform the each term of sum into the term $K_2  K_1 K_2$:             
\ba
K_2^2  K_1 + K_1 K_2^2  &=&2 K_2 K_1 K_2 + (q -q^{-1}) [K_2 ; A_+B_-  - B_+A_-]=
\nonumber \\
\;   &=&(q^{2} +q^{-2}) K_2 K_1 K_2   - (q -q^{-1})^2 (p_1 K_2 +t_1 K_1 +t_0),
\lab{17} \ea
where  $t_1 = (q +q^{-1})^2 a_2 b_1; t_0 = (q +q^{-1})(Q_B a_2 q^{-C_0}  + Q_A b_1 q^{ C_0}).$

Here, a great help in getting the final formula is to use the following expression for term $ K_2 K_1 K_2 $:
\ba
K_2  K_1 K_2 &=& [P_2 + (\pm) q^{ \Delta }] q^{ -\Delta } [P_2 + (\pm) q^{ \Delta }]  = 
\nonumber \\
\;                 &=&P_2^2  q^{-\Delta } +  P_2 q^{ -\Delta } (\pm) q^{ -\Delta } +  (\pm) P_2 +  (\pm)  (\pm) q^{ \Delta }  = 
\nonumber \\
\;                 &=&  p_1 K_2 +p_2^2 q^{ 3 \Delta } +  p_1 p_2 q^{ 2\Delta } +  
\nonumber \\
			&+&  (q^2A_+^2B_-^2+ A_+A_-B_-B_+ + A_-A_+B_+B_- + q^{-2}A_-^2B_+^2) q^{ \Delta }.
\lab{18} \ea

Thus, in the two-operator approach to the Askey-Wilson algebra \ci{la1,daha}, we obtain the following relations:
\ba
K_1^2  K_2 + K_2 K_1^2  - (q^{2} +q^{-2}) K_1 K_2 K_1+(q -q^{-1})^2 (p_2  + p_1  K_1)= 0,
\nonumber \\
\;   K_2^2  K_1 + K_1 K_2^2   - (q^{2} +q^{-2}) K_2 K_1 K_2   + (q -q^{-1})^2 (p_1 K_2 +t_1 K_1 +t_0)= 0,
\lab{myawa} \ea

Introducing the procedure of "q-mutation" for
arbitrary operators $L,M$
\be
[L,M]_{q}\equiv q LM - q^{-1} ML
\lab{qmut} \ee

we get  a special case of Askey-Wilson algebra $AW(3)$ with three generators for the equation \re{myawa}, i. e. 
 operators $K_1,K_2$ together with
their q-mutator $K_3$ obey the following algebra
\ba
[K_1,K_2]_{q}&=&K_3,
\nonumber \\
\quad[K_2,K_3]_{q}&=&BK_2+C_1K_1+D_1,
\nonumber \\
\quad[K_3,K_1]_{q}&=&BK_1+C_2K_2+D_2,
\lab{AW3} \ea
where $B,C_{1,2},D_{1,2}$ are the structure constants of the algebra \re{AW3}:
\ba
B&=& (q -q^{-1})^2 p_1=  (q -q^{-1})^2 (Q_B q^{C_0}  + Q_A q^{- C_0}),
\nonumber \\
C_1 &=& (q -q^{-1})^2 t_1  = (q^{2}-q^{-2})^2 a_2 b_1,\quad C_2=0;
\nonumber \\
D_1 &=& (q -q^{-1})^2 t_0  =(q -q^{-1})^2 (q +q^{-1})(Q_B a_2 q^{-C_0}  + Q_A b_1 q^{ C_0}),
\nonumber \\ 
D_2 &=& (q -q^{-1})^2 p_2 =(q -q^{-1})^2 (q +q^{-1}) d .
\lab{str} \ea 

The  Askey-Wilson algebra with three generators
$AW(3)$ was introduced and  studied in \ci{zh1}. The Casimir
operator $\hat Q$
commuting with all the generators $K_1,K_2,K_3$ of the our algebra
has the expression
\ba
\hat Q= \fac 12 \lbr K_3,\tilde K_3\rbr  +(q^{2} +q^{-2}) C_1K_1^2 + B \lbr K_1,K_2\rbr + (q +q^{-1})^2  (D_1K_1+D_2K_2) 
\lab{AWQ} \ea
where the symbol $\lbr.,.\rbr$ stands for the anticommutator $\{a,b\}=ab+ba$ and
$\tilde K_3$ is the
"dual" generator:
\be
\quad \tilde K_3 = [K_1,K_2]_{-q}=q^{-1}K_1K_2-qK_2K_1
\lab{dual} \ee 

Thus,  we considered algebraic treatment of hidden symmetry for the general case of addition of nonlinear algebras including $sl_q(2)\oplus sl_q(2)$  and 
two types of q-oscillator algebra. To stress this aspect as the goal of this paper the simple example are given below.  
Also we are specially pointed that for the first time many different special types of addition rule  were discussed in \ci{zh4}-\ci{zh5}.

At the beginning it follows to note that famost Jordan-Wigner realization for $SU(1,1)$ or $SU(2)$, based on the two independent, but same structure set of harmonic oscillator operators  and its q-analog, is impossible in discussed approach. So,  we can unusual add (but fully justified in given approach!) two different types of q-oscillator algebra: $eu_q^-$ and $eu_q^+$.  In symbolic form  this  non-commutative addition can be written both as
\be   
eu_q^+ \oplus eu_q^- =  sl_q(2)   \quad or \quad  (a_2, 0) \oplus (0, b_1) = (a_2, b_1)
\lab{symb} \ee

and   as
\be   
eu_q^- \oplus eu_q^+  =M(2)   \quad \ci{m2} \quad or \quad  (0, d=b_1) \oplus (d=a_2) = (0, 0).
\lab{symb} \ee

According to the first variant of the addition we can consider the algebra $ sl_q (2) $ itself as the resulting algebra from the two different types of q-Bose algebras that has its two parameters $a_2$, $b_1$. Here hidden symmetry of the Hahn problem is determinated by a special case of Askey-Wilson algebra $AW(3)$:
\ba
[K_1,K_2]_{q}&=&K_3,
\nonumber \\
\quad[K_2,K_3]_{q}&=&BK_2+C_1K_1+D_1,
\nonumber \\
\quad[K_3,K_1]_{q}&=&BK_1,
\lab{AW3a} \ea
where $B,C_{1,2},D_{1,2}$ are the structure constants of the algebra \re{AW3a}:
\ba
B&=& (q -q^{-1})^2 p_1=  (q -q^{-1})^2 (Q_B q^{C_0}  + Q_A q^{- C_0}),
\nonumber \\
C_1 &=& (q -q^{-1})^2 t_1  = (q^{2}-q^{-2})^2 a_2 b_1,\quad C_2=D_2 =0;
\nonumber \\
D_1 &=& (q -q^{-1})^2 t_0  =(q -q^{-1})^2 (q +q^{-1})(Q_B a_2 q^{-C_0}  + Q_A b_1 q^{ C_0}),
\lab{str} \ea 

According to the results of the work \ci{zh1} the overlaps between two eigenbases $\psi_p$ and $\phi_s$ as
the Clebsch-Gordan coefficientss  are expressed
in terms of special case for  the Askey-Wilson polynomials \ci{AW1}  -  q-analog of Kravchuk, Meixner, Charlier polynomials (basic hypergeometric function $_3\Phi_2$ or $_1\Phi_2$ for $C_1=0$).

The second variant of the addition are represented  the resulting algebra $M(2)$ with only one parameter $d=a_2=b_1$. Here the hidden symmetry algebra of the Hahn problem  is  the following special case of Askey-Wilson algebra $AW(3)$:
\ba
[K_1,K_2]_{q}&=&K_3,
\nonumber \\
\quad[K_2,K_3]_{q}&=&BK_2,
\nonumber \\
\quad[K_3,K_1]_{q}&=&BK_1+D_2,
\lab{AW3b} \ea
where $B,C_{1,2},D_{1,2}$ are the structure constants of the algebra \re{AW3b}:
\ba
B&=& (q -q^{-1})^2 p_1=  (q -q^{-1})^2 (Q_B q^{C_0}  + Q_A q^{- C_0}),
\nonumber \\
D_2 &=& (q -q^{-1})^2 p_2 =(q -q^{-1})^2 (q +q^{-1}) d,
\nonumber \\
C_1 &=& D_1=C_2=0,
\lab{str} \ea 
According to the results of the work \ci{zh1} the overlaps between two eigenbases $\psi_p$ and $\phi_s$ as
the Clebsch-Gordan coefficients  are expressed
in terms of special case of q-Hahn polynomials.

\vspace{6mm}
\section{The Clebsch-Gordan coefficients in the Hahn problem} 

On the base of previously reviewed hidden symmetry now let's solve the the Hahn problem or find  the corresponding  overlaps Clebsch-Gordan coefficients between two eigenbases $\psi_p$  and $\phi_s$ 
for the operatots $K_1,K_2$  respectively:
\ba
K_1\psi_p&=& \lambda_p\psi_p,  \quad K_1  = q^{- \Delta} =q^{- A_0 + B_0},
\nonumber \\
\lab{K1s} \\  
K_2\phi_s &=& \mu_s\phi_s, \quad  K_2  =  (q A_+B_-  + q^{-1} B_+A_-) q^{ \Delta } + p_2 q^{2  \Delta }  + p_1 q^{ \Delta } + p_0;
\lab{K2s} 
 \ea

Recall first that $\mathcal{C}=\mathcal{A}\oplus\mathcal{B}$ and the following relations hold
 \ba
 \lab{Actions-2}
A_0 e_{n_a}^{(\mu_a)} & = (n_a + \mu_a) e_{n_a}^{(\mu_a)}, 
\quad A_{+} e_{n_a}^{(\mu_a)} = r_{n_a+1} e_{n_a+1}^{(\mu_a)}, 
\quad A_{-} e_{n_a}^{(\mu_a)} = r_{n_a} e_{n_a-1}^{(\mu_a)},
\lab{a}\\
B_0 e_{n_b}^{(\mu_b)} & = (n_b + \mu_b) e_{n_b}^{(\mu_b)}, 
\quad B_{+} e_{n_b}^{(\mu_b)} = r_{n_b+1} e_{n_b+1}^{(\mu_b)}, 
\quad B_{-} e_{n_b}^{(\mu_b)} = r_{n_b} e_{n_b-1}^{(\mu_b)},
\lab{b}\\
C_0 e_{n_c}^{(\mu_c)} & = (n_c + \mu_c) e_{n_c}^{(\mu_c)}, 
\quad C_{+} e_{n_c}^{(\mu_c)} = r_{n_c+1} e_{n_c+1}^{(\mu_c)}, 
\quad C_{-} e_{n_c}^{(\mu_c)} =r_{n_c} e_{n_c-1}^{(\mu_c)},
 \lab{c}\ea
The first set of eigenvectors $\psi_p$ correspond to the elements of the direct product basis $\psi_p \equiv e_{n_a}^{(\mu_a)} \otimes e_{n_b}^{(\mu_b)}$.
This basis vectors of the direct product are characterized as eigenvectors of the operators
\be
\label{Label-1}
\hat Q_{A},\;\; A_{0},\;\; \hat Q_{B}, \;\; B_{0}
\ee
with eigenvalues
\be
Q(\mu_a),\;\; n_a+\mu_a, \;\; Q(\mu_b),\;\; n_b+\mu_b \;\;
\ee
respectively. 
The second set of eigenvectors $\phi_s$ is identified as should be to the coupled basis elements $e_{n_c}^{(\mu_c)}$, which are the eigenvectors of
\be
\label{Coupled-Basis}
\hat Q_{C},\;\; C_0 \equiv A_0+B_0,
\ee
with eigenvalues
\be
Q(\mu_{c}), \;\; n_c+\mu_{c} \equiv n_a+ n_b+\mu_{a} +\mu_{b},
\ee
respectively. The direct product basis is related to the coupled basis by a unitary transformation whose matrix elements are called Clebsch-Gordan coefficients. 
These overlap coefficients will be zero unless
\ba
n_{a}+n_{b} & \equiv N =n_c  + \mu_{c} - \mu_{a} - \mu_{b}.
\ea
Since $n_c$ is an integer, it follows that
\ba
\lab{CND-2}
\mu_{c} & =  \mu_{a} + \mu_{b} + x;  \quad  n_c  + x = n_{a}+n_{b} 
\ea
where $x\in\{0,\ldots,N\}$ for a given value of $N=n_{a}+n_{b}$. 

We may hence write
\be
\label{Relation}
e_{n_c}^{(\mu_c)}=\sum_{n_a, n_b}C_{\mu_a, n_a; \mu_b, n_b}^{ \mu_c, n_c}\; e_{n_a}^{(\mu_a)} \otimes e_{n_b}^{(\mu_b)},
\ee
where 
\be
\label{CG-Def}
C_{\mu_a, n_a; \mu_b, n_b}^{ \mu_c, n_c}=\left\langle {e_{n_a}^{(\mu_a)} \otimes e_{n_b}^{(\mu_b)}} | e_{n_c}^{(\mu_c)} \right\rangle  \\
\ee
are the Clebsch-Gordan coefficients of $sl_{q}(2)$.

Let's show that this Clebsch-Gordan coefficients  are expressed
by q-Hahn polynomials or basic hypergeometric function $_3\Phi_2$. Also note that the explicit expression for the Clebsch-Gordan coefficients (\ref{CG-Def}) is known \cite{zh4}-\cite{zh5}, \cite{qccc}, hence only a short derivation using a recurrence relation is presented. By definition of the coupled basis states (\ref{Coupled-Basis}), one has
\be
\label{Side-1}
Q(\mu_{c}) C_{\mu_a, n_a; \mu_b, n_b}^{ \mu_c, n_c} \equiv
[- c_1 q^{2\mu_c - 1 } - c_2 q^{1 - 2\mu_c}]C_{\mu_a, n_a; \mu_b, n_b}^{ \mu_c, n_c}=
 \left\langle {e_{n_a}^{(\mu_a)} \otimes e_{n_b}^{(\mu_b)}} | \hat Q_{C} | e_{n_c}^{(\mu_c)} \right\rangle
\ee
On the other hand, upon using (\ref{Cas-Full}) and the actions (\ref{Actions-2}), i. e. substituting the expressions (\ref{res}) for $C_{\pm}$ and $C_0$ into the right-hand side of (\ref{Side-1}), one finds
\ba
\nonumber
 \left\langle {e_{n_a}^{(\mu_a)} \otimes e_{n_b}^{(\mu_b)}} | \hat Q_{C} | e_{n_c}^{(\mu_c)} \right\rangle
 & =     \left[ p_2 q^{ 2(n_a-n_b +\mu_a-\mu_b) }  + p_1 q^{n_a-n_b +\mu_a-\mu_b} + p_0 \right] C_{\mu_a, n_a; \mu_b, n_b}^{ \mu_c, n_c}
\nonumber \\
& + r_{n_a+1}r_{n_b-1} q^{ n_a - n_b +\mu_a-\mu_b +1 }   
C_{\mu_a, n_a+1; \mu_b, n_b-1}^{ \mu_c, n_c}
\\
\label{Side-2}
& +  r_{n_a-1}r_{n_b+1} q^{ n_a -n_b +\mu_a-\mu_b -1 }   
C_{\mu_a, n_a-1; \mu_b, n_b+1}^{ \mu_c, n_c}
\ea
   
For a given value of $N=n_a+n_b$, taking $n_a=n$ and $n_b=N-n$, one can use the conditions (\ref{CND-2}) to make explicit the dependence of $\mathcal{C}$ on $x$:
$$
C_{\mu_a, n_a; \mu_b, n_b}^{ \mu_c, n_c}=\omega\,P_{n}(x;\mu_a,\mu_b;N),
$$
where $\omega=C^{ \mu_a  + \mu_b +  x; N-x}_{\mu_a, 0; \mu_b, N}$ and $P_0(x)=1$.
With these definitions, it follows from (\ref{Side-1}) and (\ref{Side-2}) that $P_{n}(x)$ satisfies the three-term recurrence relation
\ba
\lab{Recu-Norm}
 \lambda(x) P_{n}(x)&=Z_n\,P_{n}(x) +W_{n}\,P_{n-1}(x)+W_{n+1}\,P_{n+1}(x),
\ea
where 
$$\lambda(x)=(-c_2 q^{1-2\mu_a -2\mu_b}) (q^{-2x}+ c_1q^{4\mu_a +4\mu_b - 2 +2x})$$
and 
matrix elements $W_{n}$ and $Z_{n}$ can be rewritten in the form:
\ba
W_{n}^2&=&q^{2(n-N +\mu_a -\mu_b-1)} r_n^2 r_{n-N+1}^2 
\nonumber \\
&\equiv &(1-q^{2n})(1-rq^{2n})(1-q^{2(n-N-1)})(1-sq^{2(n-N)});
\nonumber \\
Z_{n}&=&(D q^{2n} + E q^{4n})
\ea
These recurrent relations with the q-dependence of  $W_{n}$ and $Z_{n}$ directly indicates on q-Hahn polynomials  \ci{AW1}-\ci{AW2}.
Omitting the details of calculation, we present the some results concerning the connection between parameters for the q-Hahn polynomial, the the additions algebras and the
$sl_q(2)$ :
\ba
\lab{corpar} 
r &=&{a_1\over a_2} q^{4\alpha-2},\qquad  s ={b_2\over b_1} q^{4\beta}, \qquad a_1 \equiv -b_2
\nonumber \\
D&=&a_2 q^{-2(\mu_a +\mu_b) }
\left[
q^{-2N}+ s^{-1} r q^2 +rq^2 (q^2+q^{-2N})
\right]
\qquad E=a_1  (q+q^{-1}) q^{2(\mu_a  - \mu_b -N)}
\nonumber \\
\ea
and
\ba
P_{n}(x) &=& h_n
\; _3\Phi_2 \Biggl( {q^{-2n},\; q^{-2x},\;  s^{-1} r q^{2x}\atop
q^{-2N}, \; r q}; q^2 \Biggl| q^2 \Biggr),
\ea
where  $h_n$ is
some normalization factor and $_3\Phi_2$ are the basic hypergeometric function (for details see \ci{zh4}, \ci{qccc}).

\section{Conclusion}
 We have shown that a special case of Askey-Wilson algebra $AW(3)$ with three generators serves as a hidden symmetry
underlying the Hahn problem for the quantum algebra $sl_q(2)$. 
On the base of this
hidden symmetry 
the  corresponding  Clebsch-Gordan coefficients 
in terms of the q-Hahn polynomials is found.
Other most properties
 of  these coefficients (symmetry, generating functions, recurrent relations)
can be automatically derived  for all the possible representation
series, not just the discrete series
$D_\mu^+$ discussed in the article.

In future publications, the authors intend to apply this result to find exactly solvable physical problems.


\newpage

\begin{thebibliography}{40}

\bibitem{zh1} A.~S. Zhedanov, ''Hidden symmetry'' of Askey-Wilson polynomials, {\it Theoret.  and Math. Phys.}, {\bf 89(2)}, 1146--1157, 1991.

\bibitem{zh2}	 L.~Vinet,  and A.~Zhedanov, A unified algebraic underpinning for the Hahn polynomials and rational functions, 
{\it LANL, Cornell University Library}, {\bf https://arxiv.org/abs/1808.09518} (accessed 14 February 2019). 

\bibitem{zh3} A.~Zhedanov, Hidden Symmetry Algebra and Overlap Coefficients for Two Ring-Shaped Potentials
{\it Journal of Physics A General Physics}  {\bf 26(18)}, 4633-4641, (1999).
https://doi.org/10.1088/0305-4470/26/18/027

\bibitem{qH} Luc Frappat, Julien Gaboriaud, Eric Ragoucy, Luc Vinet, The q-Higgs and Askey-Wilson algebras,
{\it LANL, Cornell University Library}, {\bf https://arxiv.org/abs/1808.09518} (accessed 14 February 2019).
{\it Nuclear Physics } {\bf B 944} 114632  (2019)
https://doi.org/10.1016/j.nuclphysb.2019.114632

\bibitem{zh4} Ya. I. Granovskii, A. S. Zhedanov and O. B. Grakhovskaya,
Addition rule for nonlinear algebras
{\it Phys.Lett.}, {\bf B278}, 85-88 (1992).

\bibitem{zh5} Alexei Zhedanov,
Q rotations and other Q transformations as unitary nonlinear automorphisms of quantum algebras
{\it Journal of Mathematical Physics}  {\bf 34(6) }, 2631 (1993).
https://doi.org/10.1063/1.530088

\bibitem{maj} S. Majid, Quasitriangular Hopf algebras and Yang-Baxter equations
{\it Int.J.Mod.Phys.}, {\bf A5(1)}, 1-91, (1990).


\bibitem{la1} A. Lavrenov, 
On Askey--Wilson algebra,
in Quantum Groups and Integrable Systems, II (Prague, 1997),
{\it Czechoslovak J. Phys.}, {\bf 47(12)}, 1213-1219, (1997).
http://dx.doi.org/10.1023/A:1022821531517


\bibitem{daha} Tom H. Koornwinder, Marta Mazzocco, Dualities in the q-Askey scheme and degenerate, 
{\it LANL, Cornell University Library}, {\bf https://arxiv.org/abs/1803.02775} (accessed 25 September 2018).
https://doi.org/10.1111/sapm.12229
 {\it Studies in Applied Mathematics}, {\bf 141(4)}, 424-473, ( 2018).


\bibitem{m2} E. G. Kalnins, W. Miller, Jr., and S. Mukherjee, 
Models of Q-Algebra Representations: The Group of Plane Motions
{\it SIAM J.Math.Anal.}, {\bf 25(2)}, 513–527 (1994) 
https://doi.org/10.1137/S0036141092224613


\bibitem{qccc} Vincent X. Genest, Plamen Iliev, and Luc ~Vinet,
Coupling coefficients of suq(1,1)
and multivariate q-Racah polynomials,
{\it LANL, Cornell University Library}, {\bf https://arxiv.org/abs/1702.04626} (accessed 14 February 2019). 
{\it Nuclear Physics B}, {\bf 927}, 97-123 (2018) 
https://doi.org/10.1016/j.nuclphysb.2017.12.009



\bibitem{AW1} R.Askey and J.Wilson, A Set of Orthogonal Polynomials That Generalize the Racah Coefficients or $6 - j$ Symbols {\it SIAM J.Math.Anal.}, {\bf 10(5)},
1008–1016 (1979); 
https://doi.org/10.1137/0510092


\bibitem{AW2} R.Askey and J.Wilson, Some basic hypergeometric orthogonal polynomials that generalize Jacobi polynomials {\it Mem.Am.Math.Soc.}, {\bf 54 (319)}, 1-55 (1985).
DOI: http://dx.doi.org/10.1090/memo/0319



\end{thebibliography}
\end{document}